\documentclass{amsart}
\usepackage{amsmath}
\usepackage{amssymb}
\usepackage{amsthm}
\usepackage{dsfont}
\usepackage{lineno}
\usepackage{subfigure}
\usepackage{graphicx}
\usepackage{multirow}
\usepackage{color}
\usepackage{xspace}
\usepackage{pdfsync}
\usepackage{hyperref} 
\usepackage{algorithmicx}
\usepackage{algpseudocode}
\usepackage{algorithm}
\usepackage{mathtools}
\usepackage{enumitem}
\graphicspath{{Figures/}}

\DeclareMathOperator*{\argmin}{argmin}

\DeclarePairedDelimiter\floor{\lfloor}{\rfloor}
\newcommand{\bq}{\begin{equation}}
\newcommand{\eq}{\end{equation}}
\newcommand{\R}{\mathbb{R}}

\newcommand{\abs}[1]{\left\vert#1\right\vert}
\newcommand{\norm}[1]{\left\Vert#1\right\Vert}

\newcommand{\G}{\mathcal{G}}
\newcommand{\bO}{\mathcal{O}}

\newcommand{\Dt}{\mathcal{D}}

\newcommand{\Sf}{\mathbb{S}^{2}}

\newcommand{\Zf}{\mathcal{Z}}

\newcommand{\MA}{Monge-Amp\`ere\xspace}

\newcommand*{\avint}{\mathop{\ooalign{$\int$\cr$-$}}}

\algnewcommand{\LineComment}[1]{\State \(\triangleright\) #1}

\newtheorem{theorem}{Theorem}

\theoremstyle{lemma}

\newtheorem{remark}[theorem]{Remark}

\theoremstyle{remark}

\makeatletter
\newcommand\appendix@section[1]{%
\refstepcounter{section}%
\orig@section*{Appendix \@Alph\c@section: #1}%
}
\let\orig@section\section
\g@addto@macro\appendix{\let\section\appendix@section}
\makeatother

\begin{document}

\title[Reflector Antenna Problem via Optimal Transport]{A Convergent Numerical Method for the Reflector Antenna Problem via Optimal Transport on the Sphere}

\author{Brittany Froese Hamfeldt}
\address{Department of Mathematical Sciences, New Jersey Institute of Technology, University Heights, Newark, NJ 07102}
\email{bdfroese@njit.edu}
\author{Axel G. R. Turnquist}
\address{Department of Mathematical Sciences, New Jersey Institute of Technology, University Heights, Newark, NJ 07102}
\email{agt6@njit.edu}

\thanks{The first author was partially supported by NSF DMS-1619807 and NSF DMS-1751996. The second author was partially supported by  an NSF GRFP}

\begin{abstract}
We consider a PDE approach to numerically solving the reflector antenna problem by solving an Optimal Transport problem on the unit sphere with cost function $c(x,y) = -2\log \left\Vert x - y \right\Vert$. At each point on the sphere, we replace the surface PDE with a generalized Monge-Amp\`ere type equation posed on the local tangent plane. We then utilize a provably convergent finite difference scheme to approximate the solution and construct the reflector. The method is easily adapted to take into account highly nonsmooth data and solutions, which makes it particularly well adapted to real-world optics problems.  Computational examples demonstrate the success of this method in computing reflectors for a range of challenging problems including discontinuous intensities and intensities supported on complicated geoemtries.
\end{abstract}

\date{\today}    
\maketitle
\section{Introduction}\label{sec:intro}
Advances in light emitting diode (LED) technology in recent years have allowed for more flexibility in the engineering of freeform lenses using plastics in light illumination problems. In this article, we focus on the reflector antenna problem, which involves designing a reflector to reshape a point source onto a prescribed output in the far-field.  
On the theoretical side, a major advance in understanding freeform geometric optics problems has been gained by reformulating the problem as a fully nonlinear partial differential equation (PDE) of Monge-Amp\`{e}re type. In the particular case of the reflector antenna problem, this PDE is posed on the sphere.  The curved geometry, nonlinearity of the equation, and singular terms within the PDE make this a challenging problem to solve numerically.  

In this article, we propose a new method for the design of the reflector surface that relies on recent advances by the authors in the numerical approximation and analysis of \MA type equations on the sphere~\cite{HT_OTonSphere2}.  We emphasize that this new method comes with theoretical guarantees of convergence, even in settings involving very non-smooth output intensities~\cite{HT_OTonSphere}.

Computational approaches to solving optical design problems can be roughly divided into three basic categories: (1) techniques that use a ray-mapping to design the optical surface, (2) methods that approximate the optical surfaces by supporting quadrics, and (3) methods that represent the optical surface through the solution to an optimal transportation problem.

The ray-mapping approach generally involves a two-step procedure.  In the first step, a ray mapping is produced between the input and output intensities. In the second step, the laws of reflection and/or refraction are employed to construct a surface that achieves this ray mapping as nearly as possible. Several methods based on this general approach are available including~\cite{Bruneton_lens,Desnijder_raymapping,FFL_optics,Fournier_reflector,Parkyn_lenses}.  A downside to this general approach is that it can be difficult to theoretically justify the existence of an optical surface that exactly produces the desired ray mapping.

Oliker's method of supporting quadrics involves representing the optical surface via supporting ellipsoids or hyperboloids~\cite{Oliker_nearfield,Oliker_SQM}. The simple optical properties of these quadrics is used to produce a pixelated version of the desired target.  This approach has the advantage of being theoretically well-founded, but can be costly to implement in practice.

Finally, the solution to many optical design problems can be obtained directly through the solution of a corresponding optimal transportation problem.  That is, if $f_1$ represents the input intensity and $f_2$ the desired output intensity, it is necessary to solve a problem of the form
\bq\label{eq:OT}
\min\limits_{T_\# f_1 = f_2} \int_{\text{supp}(f_1)} c(x,T(x)) f_1(x) dx.
\eq
where $c(x,y)$ is the cost of transporting a unit of mass from $x$ to $y$ and $T_\# f_1 = f_2$ indicates that
\bq\label{eq:massCons} \int_A f_1(x)\,dS(x) = \int_{T(A)} f_2(y)\,dS(y) \eq
for every measurable $A \subset \Sf$.

Many optical inverse problems have yielded fruitful interpretations via optimal transport by deriving an appropriate cost function $c(x,y)$~\cite{YadavThesis}. To give a simple example, a parallel-in, far-field out setup yields the cost function $c(x,y) = \frac{1}{2}\left\Vert x - y \right\Vert^2$, where $x,y \in \mathbb{R}^2$.  The reflector antenna problem considered in this article
 has a slightly more challenging set-up in that the cost function $c(x,y) = -2\log \left\Vert x - y \right\Vert$ is unbounded and the intensity functions $f_1, f_2$ are supported on $\mathbb{S}^2$ (the unit $2$-sphere), as opposed to subsets of Euclidean space~\cite{GangboOliker,OlikerNewman,Wang_Reflector,Wang_Reflector2}. 

One approach to solving optimal transport problems in optical design is to use optimization techniques, including linear assignment~\cite{Doskolovich_farfield} and linear programming~\cite{GlimmOliker_SingleReflector}.  This approach has the advantage of being theoretically well-understood.  However, the optimization problems typically involve a very large number of constraints and the resulting methods are computationally complex. 

In many cases, the solution to the optimal transport problem can also be obtained through the solution of a fully nonlinear partial differential equation of \MA type, which has the general form
\bq\label{eq:MA}
\det(D_{xx}^2(u(x)+A(x,\nabla u(x)))) = H(x,\nabla u(x))
\eq
subject to the constraint that
\bq\label{eq:cconvex}
D_{xx}^2(u(x)+A(x,\nabla u(x))) \geq 0,
\eq
where $M\geq 0$ means that $M$ is positive semi-definite.  In the case of a point source lens or reflector design problem, this PDE is posed on the unit sphere $\Sf$.

Recently, several methods have been proposed for solving optical design problems involving a point source via the solution of a \MA type equation.  These methods replace the PDE on the sphere with a corresponding equation on the plane by representing subsets of the unit sphere using spherical coordinates~\cite{Wu_lensdesign}, a vertical projection of coordinates onto the plane~\cite{Brix_MAOptics}, or stereographic projection~\cite{RomijnSphere}.  As the numerical solution of these \MA type equations is a very new field, many of the numerical methods used in optical design problems are not yet equipped with theoretical guarantees of convergence.

In the present article, the solution to the reflector antenna problem is obtained by solving a \MA type equation directly on the sphere.  This has the advantage of allowing for intensity distributions supported on complicated subsets of the sphere or even the entire sphere.  Moreover, the approach is intrinsic and thus the solution to the problem will not depend on such details as the choice of the north pole.  Finally, the numerical method we use is theoretically well-justified and can be proven to converge to the correct solution of the \MA equation in a wide variety of challenging settings~\cite{HT_OTonSphere, HT_OTonSphere2}.

\section{Mathematical Approach}\label{sec:background}

Here we briefly summarize the derivation of the reflector antenna problem and its connection to optimal transport on the sphere, which leads to an equation of \MA type that can be solved using techniques from numerical PDEs. 

\begin{figure}[htp]
\includegraphics[height=8cm]{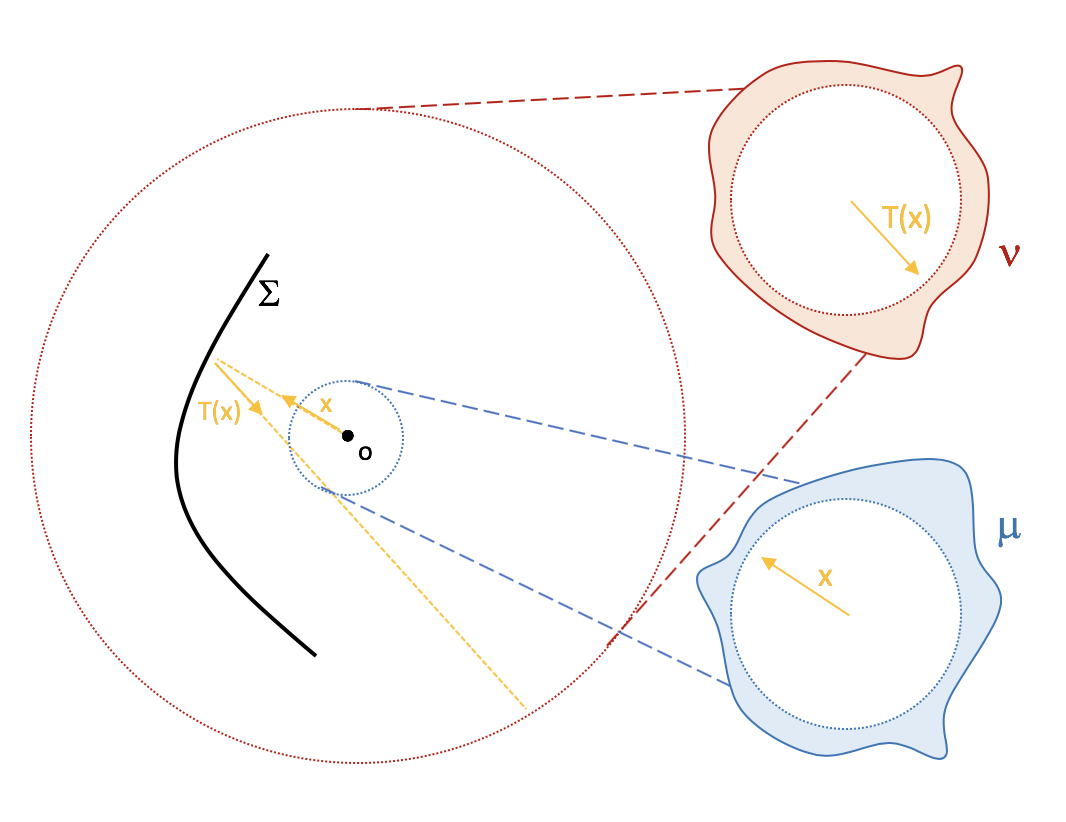}
\caption{Reflector antenna with source/detector $\mu$, reflector $\Sigma$ and target far-field intensity $\nu$. The directional vectors $x$ and $T(x)$ are unit vectors.}\label{fig:reflectorantenna}
\end{figure}

We begin by following the physical derivation in~\cite{Wang_Reflector, Wang_Reflector2}. We start with a light source or detector $\mu$ located at the origin, which is a probability measure indicating directional intensity and is supported on a set $\Omega \subset \mathbb{S}^2$. Next we consider a reflector surface $\Sigma$, which is a radial graph over the domain $\Omega$ and can be represented as
\begin{equation}
\Sigma = \left\{ x\rho(x) \mid x \in \Omega, \quad \rho>0 \right\}
\end{equation}
where $\rho: \Omega \rightarrow \mathbb{R}$ is a non-negative function indicating the distance between the reflector surface and the origin. The light from the source $\mu$ in the direction $x$ bounces off the reflector $\Sigma$ without any refraction or absorption and travels in the direction $T$ following the law of reflection. Over all directions this produces the far-field intensity $\nu$, which is also a probability measure indicating directional intensity and is supported on some target domain $\Omega^{*} \subset \Sf$. See Figure~\ref{fig:reflectorantenna} for a schematic of the setup.

The reflector antenna problem is thus: given source and target intensity probability distributions $\mu$ and $\nu$, respectively, find the shape of the reflector $\Sigma$ that transmits the light from the source to the target while satisfying conservation of energy.
We make the assumption that the probability densities $\mu$ and $\nu$ have density functions $f_1$ and $f_2$ respectively (so that $d\mu(x) = f_1(x)dS(x), d\nu(y) = f_2(y)dS(y)$).  Now we seek a PDE that will allow us to determine the reflector height function $\rho(x)$, which fully determines the reflector surface, in terms of the prescribed intensity functions $f_1$ and $f_2$.

\begin{figure}[htp]
\includegraphics[height=4.5cm]{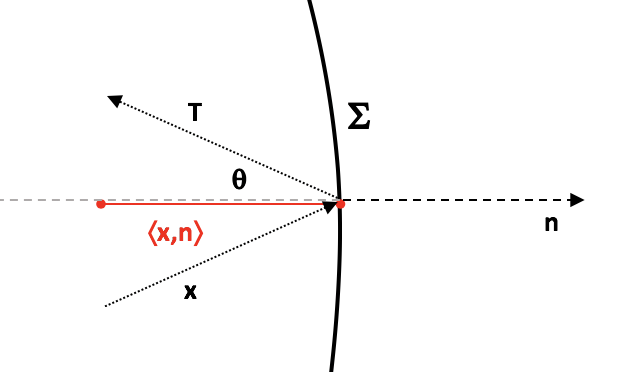}
\caption{Incident light direction $x$, reflector $\Sigma$, outward normal $n$, and outward light ray $T$}\label{fig:reflection}
\end{figure}

The first of the two physical laws that will be used to derive the governing PDE for this setup is the well known geometric law of reflection, which yields the optical map
\begin{equation}\label{eq:mapping}
T(x) = x - 2\left\langle x, n(x) \right\rangle n(x)
\end{equation}
where $n(x)$ is the outward normal to $\Sigma$ at the point $z = x \rho(x)$, $x \in \Omega$.  See Figure~\ref{fig:reflection}.  We emphasize here that this is the geometric optics limit.

The second physical law that completes the problem is the law of conservation of energy:
\begin{equation}\label{eq:conservation}
\int_{{T}^{-1}(E)} f_1(x)\,dS(x) = \int_{E} f_2(y)\,dS(y).
\end{equation}
for any Borel set $E \subset \Omega^{*}$.

By introducing local coordinates on the sphere, Wang~\cite{Wang_Reflector} observes that the unit normal $n$ can be given by
\bq\label{eq:n}
n(x) = \frac{\nabla\rho(x) - x\rho(x)}{\sqrt{\rho(x)^2+\norm{\nabla\rho(x)}^2}}.
\eq
Then the law of reflection~\eqref{eq:mapping} yields the mapping
\bq\label{eq:mapping2}
T(x) = \frac{2\rho(x)\nabla\rho(x)+\left(-\rho(x)^2+\norm{\nabla\rho(x)}^2\right)x}{\rho(x)^2+\norm{\nabla\rho(x)}^2}.
\eq
Applying the change of variables formula to the conservation of energy constraint~\eqref{eq:conservation} produces an equation of the form
\bq\label{eq:detEquation}
\det(\nabla T(x)) = f_1(x)/f_2(T(x)).
\eq
Combining these conditions yields the PDE
\begin{equation}\label{eq:opticsPDE}
\eta^{-2} \det\left( -\nabla_{i} \nabla_{j} \rho + 2 \rho^{-1} \nabla_{i} \rho \nabla_{j} \rho + (\rho - \eta)\delta_{ij} \right) = f_1(x)/f_2(T(x))
\end{equation}
where $\eta = \left( \left\vert \nabla \rho \right\vert^2 + \rho^2 \right)/2\rho$ and $\delta_{ij}$ is the usual Kronecker delta. We recognize this PDE as an equation of Monge-Amp\`{e}re type, with the usual second boundary value condition~\cite{Urbas}
\begin{equation}
T(\Omega) = \Omega^{*}.
\end{equation}

Unfortunately, there are few direct results in the literature that answer the kind of questions of existence and regularity results that are needed to design a convergent numerical method for~\eqref{eq:opticsPDE}. Instead, we extract a problem with more structure via the change of variables
\begin{equation}\label{eq:cov}
\rho = e^{-u}.
\end{equation}
Wang~\cite{Wang_Reflector2} shows that under an equivalent change of variables (modulo a sign change), the function $u$ solves the dual formulation of the optimal transport problem with cost function $\tilde{c}(x,y) = -\log(1-x\cdot y)$.

As an alternative approach, we notice that under this change of variables, the optical mapping~\eqref{eq:mapping2} becomes
\bq\label{eq:mappingu}
T(x) = \frac{-2\nabla u(x) + \left(\norm{\nabla u(x)}^2-1\right)x}{\norm{\nabla u(x)}^2+1}.
\eq
As in~\cite{HT_OTonSphere}, we regard this mapping as a function of the two variables $(x,\nabla u(x))$ and recognize it as a solution of the system
\bq\label{eq:optimalMap}
\begin{cases}
\nabla_x c(x,T(x,p)) = -p, & x \in \Sf\\
T(x,p) \in \Sf
\end{cases}
\eq
with the cost function 
%
%
%
\bq\label{eq:logCost}
c(x,y) = -2\log \norm{ x-y }.
\eq
This is precisely the optimality condition for the optimal transport problem on the sphere~\cite{Loeper_OTonSphere}.  Combined with the conservation of energy condition~\eqref{eq:conservation}, we can conclude that the optical mapping $T(x)$ is a solution of the optimal transport problem~\eqref{eq:OT}-\eqref{eq:massCons} with cost~\eqref{eq:logCost}.  Moreover, this interpretation opens up many existence, regularity, and numerical approximation results that can be used in determining the reflector surface $\Sigma$.

Loeper has studied this problem in detail~\cite{Loeper_OTonSphere}. Under mild conditions on the intensity distributions $f_1$ and $f_2$, the function $u$ (which fully determines the reflector surface) can be uniquely obtained as the solution of the following \MA type equation.
\bq\label{eq:OTPDE}
\begin{cases}
\det(D^2u + A(x,\nabla u)) = H(x,\nabla u), & x \in \Sf\\
D^2u + A(x,\nabla u) \geq 0.
\end{cases}
\eq
Here
\bq\label{eq:PDETerms}
\begin{split}
A(x,p) &= D_{xx}^2c \left( x,T(x,p) \right)\\
H(x,p) &= \abs{\det{D_{xy}^2c \left( x,T(x,p) \right)}}f_1(x)/f_2 \left( T(x,p) \right).
\end{split}
\eq
and the statement $M \geq 0$ means that $M$ is positive semi-definite.  This constraint (related to the so-called $c$-convexity of the optimal map $T$) is needed to ensure that the PDE has a unique solution (up to additive constants) and that this solution corresponds to the desired optical mapping $T$. 

We remark that the above equation describes a nonlinear relationship between the surface gradient and Hessian on the sphere.  In light of our goal of solving this equation numerically, perhaps the most challenging term is the mixed Hessian $D_{xy}^2c(x,y)$, which involves derivatives with respect to two different variables located at different points on the sphere.  However, following the derivation in~\cite{HT_OTonSphere2}, we can obtain a very simple explicit expression for this term by interpreting it as a change of area formula:
\begin{equation}\label{eq:mixedHessian}
\abs{\det{D_{xy}^2c \left( x,T(x,p) \right)}} = \frac{\left(\left\Vert p \right\Vert^2 + 1 \right)^2}{4}.
\end{equation}


A second challenge associated with the nonlinear \MA type equation~\eqref{eq:OTPDE} is that it requires the enforcement of an additional constraint that $D^2u + A(x,\nabla u) \geq 0$, which makes it difficult to directly apply standard techniques for approximating PDEs.  However, we succeed at absorbing this constraint into the PDE itself by relying on the following characterization of a positive semi-definite $n\times n$ matrix $M$~\cite{FO_MATheory}:
\bq\label{eq:matrix}
\begin{split}
\det(M) &= \min\limits_{\nu_i^T\nu_k=\delta_{ik}}\prod\limits_{j=1}^n \nu_j^TM\nu_j \\
  &= \min\limits_{\nu_i^T\nu_k=\delta_{ik}}\prod\limits_{j=1}^n \max\{\nu_j^TM\nu_j, 0\}.
\end{split}
\eq
Here $\delta_{ij}$ denotes the Kronecker delta function and this involves a minimization over all orthogonal coordinate frames for $\R^n$.  By observing that $\nu_j^TM\nu_j \geq 0$ for any positive semi-definite matrix $M$, we can include this condition directly in the operator instead of requiring it to be specified as a separate constraint.  This allows us to reformulate the system~\eqref{eq:OTPDE}-\eqref{eq:PDETerms} as the following unconstrained PDE.
\bq\label{eq:PDEReformulated}
\begin{split}
F(x,\nabla u(x), D^2u(x)) &\equiv
\min\limits_{\nu_1\cdot\nu_2 = 0}\prod\limits_{j=1}^2 \left. \max\left\{\frac{\partial^2(u(x)-2\log\norm{x-y})}{\partial\nu_j^2},0\right\} \right|_{y = T(x,\nabla u(x))}\\ &\phantom{=}- \frac{\left(\norm{\nabla u(x)}^2+1\right)^2 f_1(x)}{4f_2(T(x,\nabla u(x))}\\
 &= 0.
\end{split}
\eq

\section{Numerical Method}\label{sec:method}
We now describe the algorithm we use to construct the reflector surface $\Sigma$.  The algorithm hinges on the numerical solution of the nonlinear PDE~\eqref{eq:PDEReformulated}.  For fully nonlinear PDEs, it is well known that consistent and stable numerical methods may nevertheless fail to compute the correct solution. In fact, because the function $u$ is unique only up to additive constants, even fairly sophisticated numerical methods can fail to find any solution at all.  The method we describe here is inspired by a numerical scheme recently designed by the authors, which is equipped with a proof of convergence to the physically meaningful solution of the optimal transport problem.  We summarize the scheme here, and refer to~\cite{HT_OTonSphere,HT_OTonSphere2} for complete details and analysis.

\subsection{Algorithm}
We begin with a high-level overview of the algorithm.  Details will be expanded on in the following subsections.

Our starting point is a finite set of $N$ grid points $\G \subset \Sf$ that discretize the unit sphere, and the intensity distributions $f_1$ and $f_2$ that are supported on domains $\Omega\subset\Sf$ and $\Omega^*\subset\Sf$ respectively.  We let $d_{\Sf}(x,y)$ denote the usual geodesic distance between points $x,y$ on the sphere.

To the grid $\G$, we associate a number $h$ that indicates the overall spacing of grid points.  More precisely,
\begin{equation}\label{eq:h}
h = \sup\limits_{x\in\Sf}\min\limits_{y\in\G^h} d_{\Sf}(x,y) = \bO\left(N^{-1/2}\right).
\end{equation}
In particular, this guarantees that any ball of radius $h$ on the sphere will contain at least one discretization point.

Now we seek a finite difference approximation of the form
\bq\label{eq:fd}
F^h(x,u;f_1,f_2) = 0, \quad x \in \G
\eq
that approximates the original PDE~\eqref{eq:PDEReformulated}.  Our goal is to construct an approximation with the properties that (1) a solution $u^h$ exists and (2) the solution is close to the solution $u$ of the original PDE.  Our earlier work~\cite{HT_OTonSphere,HT_OTonSphere2} provides a framework for doing this.  In the most challenging settings, this requires some initial preprocessing of the data $f_1, f_2$, but then provides us with an algorithm that is guaranteed to produce a reflector surface $\Sigma^h$ that is close to the desired reflector $\Sigma$.
See Algorithm~\ref{alg:reflector}.

\begin{algorithm}[h]
\caption{Computing the reflector surface $\Sigma$}
\label{alg:reflector}
\begin{algorithmic}[1]
\State Preprocess data \[f_2^\epsilon \leftarrow \text{Regularize}(f_2).\]
\State Iterate
\[ u^h_{n+1} = u^h_n + k \left(F^h(x,u^h_n;f_1,{f_2^\epsilon}) - \sqrt{h}u^h_n(x) \right) \]
to steady state.
\State Normalize solution
\[ u^h(x) \leftarrow u^h(x) - \avint_{\phantom{==}\Sf}u^h(x)\,dS(x). \]
\State Construct reflector
\[ \Sigma^h = \left\{xe^{-u^h(x)} \mid x \in \Omega \cap \G\right\}. \]
\end{algorithmic}
\end{algorithm}

\subsection{Discretization}
We now consider a fixed grid point $x_0 \in \G$ and a grid functions $u:\G\to\R$ and explain how we obtain the value of $F^h(x_0,u;f_1,f_2)$; we refer to~\cite{HT_OTonSphere2} for further details.

We begin by projecting grid points close to $x_0$ onto the tangent plane at $x_0$.  That is, we consider the set of relevant discretization points
\bq\label{eq:neighbourCandidates}
\Zf(x_0) = \left\{z = \text{Proj}(x;x_0) \mid x \in \G, d_{\Sf}(x,x_0) \leq \sqrt{h}\right\}.
\eq
The projection is accomplished using geodesic normal coordinates, which are chosen to preserve the distance from $x_0$ (i.e. $d_{\Sf}(x,x_0) = \|x_0 -  \text{Proj}(x;x_0)\|$).  This prevents any distortions that would affect the second order terms in the PDE~\eqref{eq:PDEReformulated}.
\bq\label{eq:projection}
\text{Proj}(x;x_0) = x_0\left(1-d_{\Sf}(x_0,x)\cot d_{\Sf}(x_0,x)\right) + x \left(d_{\Sf}(x_0,x)\csc d_{\Sf}(x_0,x)\right).
\eq

The form of~\eqref{eq:PDEReformulated} indicates that we will need to approximate derivatives along various directions $\nu$.  We will consider the following finite set of possible directions,
\bq\label{eq:directions}
V = \left\{\left\{(\cos(jd\theta),\sin(jd\theta)),(-\sin(jd\theta),\cos(jd\theta))\right\} \mid j=1,\ldots,\frac{\pi}{2d\theta}\right\},
\eq
where the angular resolution $d\theta = \dfrac{\pi}{2\floor{\pi/(2\sqrt{h})}}$.

For each $\nu\in V$, we need to select four grid points $x_j\in\Zf(x_0)$, $j=1, \ldots, 4$, which will be used to construct the directional derivatives in this direction.  To accomplish this, we let $\nu^\perp$ be a unit vector orthogonal to $\nu$ and represent points in $x\in\Zf(x_0)$ using (rotated) polar coordinates $(r,\theta)$ centred at $x_0$ via
\[ x = x_0 + r(\nu\cos\theta + \nu^\perp\sin\theta), x \in \Zf(x_0). \]
Then we select four points, each in a different quadrant ($Q_1, \ldots, Q_4$), that are well-aligned with the direction of $\nu$ via
\bq\label{eq:neighbours}
x_j \in \argmin\limits_{x\in\Zf(x_0)}\left\{\abs{\sin\theta} \mid \abs{\sin\theta} \geq d\theta, r \geq \sqrt{h}-2h, x\in Q_j\right\}
\eq
where $\cos\theta \geq 0$ for points in $Q_1$ or $Q_4$ and $\sin\theta \geq 0$ for points in $Q_1$ or $Q_2$.  

From here, we construct approximations of second directional derivatives (and first directional derivatives for the usual coordinate directions $(1,0)$ and $(0,1)$) of the form
\bq\label{eq:dnu}
\begin{split}
\Dt_{\nu\nu}u(x_0) &= \sum\limits_{j=1}^4 a_j(u(x_j)-u(x_0)) \approx \frac{\partial^2u(x_0)}{\partial\nu^2}\\
\Dt_{\nu}u(x_0) &= \sum\limits_{j=1}^4 b_j(u(x_j)-u(x_0)) \approx \frac{\partial u(x_0)}{\partial\nu}.
\end{split}
\eq
The coefficients in these finite difference approximations are given explicitly by
\bq\label{eq:coeffs}
\begin{split}
a_{1} &= \frac{2\sin\theta_4(\cos\theta_3\sin\theta_2-\cos\theta_2\sin\theta_3)}{r_1\det(A)}\\
a_{2} &= \frac{2\sin\theta_3(\cos\theta_1\sin\theta_4-\cos\theta_4\sin\theta_1)}{r_2\det(A)}\\
a_{3} &= \frac{-2\sin\theta_2(\cos\theta_1\sin\theta_4-\cos\theta_4\sin\theta_1)}{r_3\det(A)}\\
a_{4} &= \frac{-2\sin\theta_1(\cos\theta_3\sin\theta_2-\cos\theta_2\sin\theta_3)}{r_4\det(A)}\\
b_{1} &= \frac{\sin\theta_4 (r_2\sin\theta_3 \cos^2\theta_2 - r_3\sin\theta_2 \cos^2\theta_3)}{r_1\det(A)} \\
b_{2} &= -\frac{\sin\theta_3 (r_1\sin\theta_4 \cos^2\theta_1 - r_4\sin\theta_1 \cos^2\theta_4)}{r_2\det(A)} \\
b_{3} &= \frac{\sin\theta_2 (r_1\sin\theta_4 \cos^2\theta_1 - r_4\sin\theta_1 \cos^2\theta_4)}{r_3\det(A)} \\
b_{4} &= -\frac{\sin\theta_1 (r_2\sin\theta_3 \cos^2\theta_2 - r_3\sin\theta_2 \cos^2\theta_3)}{r_4\det(A)}
\end{split}
\eq
where 
\bq\label{eq:detA}\begin{split}\det(A) = &(\cos\theta_3\sin\theta_2-\cos\theta_2\sin\theta_3)(r_1\cos^2\theta_1\sin\theta_4-r_4\cos^2\theta_4\sin\theta_1)\\&-(\cos\theta_1\sin\theta_4-\cos\theta_4\sin\theta_1)(r_3\cos^2\theta_3\sin\theta_2-r_2\cos^2\theta_2\sin\theta_3).\end{split}\eq

Equation~\eqref{eq:PDEReformulated} contains several functions of the gradient.  We introduce the shorthand notation
\bq\label{eq:gradFuns}
g_1(p;\nu) = \left. -2\Dt_{\nu\nu}\log\norm{x_0-y}\right|_{y=T(x_0,p)}, \quad g_2(p) = \frac{\left(\norm{p}^2+1\right)^2}{4f_2(T(x_0,p))},
\eq
denote by $L_g$ the Lipschitz constant of the function $g$, and for each function define the small parameter
\bq\label{eq:epsilon}
\epsilon_g = L_g\max\limits_{j=1,\ldots,4}\frac{\abs{b_j}}{\abs{a_j}} = \bO(\sqrt{h}).
\eq
Then all functions of the gradient can be discretized using a Laplacian regularization via
\bq\label{eq:discGrad}
g^{\pm}\left(\nabla^h u(x_0)\right) = g\left(\Dt_{(1,0)}u(x_0),\Dt_{(0,1)}u(x_0)\right) \mp \epsilon_g\left(\Dt_{(1,0),(1,0)}u(x_0) + \Dt_{(0,1),(0,1)}u(x_0)\right).
\eq

This regularization allows for the construction of a monotone scheme, which is necessary for the convergence theorem in~\cite{HT_OTonSphere}. Finally, we can combine these different operators to obtain the approximation
\bq\label{eq:approx}
\begin{split}
F^h&(x_0,u;f_1,f_2) = \\ &\min\limits_{\{\nu_1,\nu_2\}\in V} \prod\limits_{j=1}^2 \max\left\{\Dt_{\nu_j\nu_j}u(x_0) + g^-_{1,\nu_j}(\nabla^h u(x_0)), 0\right\} - f_1(x_0)g_2^+\left(\nabla^h u(x_0)\right).
\end{split}
\eq

\begin{remark}
The method of~\cite{HT_OTonSphere2} in principal involves solving a problem with this approximation, verifying that the solution satisfies required Lipschitz bounds, then if necessary solving a second discrete problem to enforce the Lipschitz condition.  However, we have never seen the verification step fail in practice, and hence never actually need to solve a second discrete system. 
\end{remark}

\subsection{Computational Complexity}
Let $N$ be the total number of grid points.  
At each point $x_0\in\G$, evaluating the operator $F^h$ involves computing a minimum over the $\bO\left(1/d\theta\right) = \bO\left(1/\sqrt{h}\right)=\bO\left(N^{1/4}\right)$ pairs of vectors in $V$.

Each pair of vectors $\{\nu_1,\nu_2\}\in V$ requires the construction of two finite difference operators of the form $\Dt_{\nu\nu}$.  Computing each of these requires identifying the four neighbors $x_1, x_2, x_3, x_4$ in the stencil.

We note that selecting each of these neighboring points $x_j$ as in~\eqref{eq:neighbours} involves searching a region whose area scales like $\bO(h^2)$.  From the definition of $h$, this is guaranteed to contain at least one point, and expected to contain $\bO(1)$ points total. Thus identification of these four neighboring points can be done in $\bO(1)$ time.

Thus, given a grid function $u$, the total computational cost of evaluating the operator $F^h$ at all points in the grid is $\bO\left(N^{5/4}\right)$.

\subsection{Preprocessing of data}\label{preprocessing}
Stability and convergence of the numerical method requires at least one of the densities (denoted by $f_2$) to be strictly positive.  This is easily accomplished by choosing $\epsilon>0$ and letting
\bq\label{eq:f2pos}
\tilde{f}_2^\epsilon = (1-\epsilon)f_2 + \frac{\epsilon}{4\pi}.
\eq
As $\epsilon\to0$, the mapping of the regularized optimal transport problem converges in measure to the solution of the given problem~\cite{Villani1}, and thus we recover the desired reflector surface.

The numerical method further requires this density function to be smoothed in order to have a (discrete) Lipschitz constant that is at most $\bO\left(h^{-1/4}\right)$.  We accomplish this via a short-time evolution of the heat equation.  That is, we solve
\bq\label{eq:heat}
\begin{cases}
v_t(x,t) = \Delta v(x,t), &(x,t)\in\Sf\times(0,\sqrt{h}]\\
v(x,0) = \tilde{f}_2^\epsilon(x), &x\in\Sf
\end{cases}
\eq
where $\Delta$ is the Laplace-Beltrami operator.  We then set
\bq\label{eq:f2reg}
f_2^\epsilon(x) = v(x,\sqrt{h}).
\eq

The Laplace-Beltrami operator can be discretized using the finite difference schemes~\eqref{eq:dnu} as
\bq\label{eq:Laplace}
\Delta^h = \Dt_{(1,0),(1,0)} + \Dt_{(0,1),(0,1)}
\eq
and evolved using forward Euler
\bq\label{eq:heatDisc}
v^{n+1} = v^n + k \Delta^hv^n.
\eq
The wide stencil nature of the finite difference stencils ($\norm{x_j-x_0} = \bO(\sqrt{h})$) means that this is stable for a time step $k \leq 1/\sum\limits_j a_j = \bO(h)$.  Thus a total of $\bO\left(h^{-1/2}\right)$ time steps are needed, which leads to an overall cost of $\bO\left(N^{5/4}\right)$ that is similar to the cost of discretization.

This regularization procedure can also be applied to unbounded densities, but requires evolving the heat equation to a stopping time of $t=h^{1/6}$ to achieve the required Lipschitz bound.

\subsection{Parabolic solvers}
After discretization, we are left with the task of solving the nonlinear algebraic system
\bq\label{eq:system}
F^h(x,u;f,g) = 0, \quad x \in \G.
\eq
Here, we use an explicit  parabolic scheme of the form
\bq\label{eq:parabolic}
u_{n+1}^h(x) = u_n^h(x) + k F^h(x,u^h_n;f,g).
\eq

As discussed in~\cite{ObermanSINUM}, we can require the time step $k$ to satisfy a nonlinear CFL condition in order to guarantee convergence. In particular, choosing $k < 1/L_{F^{h}} = \bO(h^{-2})$ is sufficient, where $L_{F^{h}}$ is the Lipschitz constant of $F^h$ with respect to the arguments $u^h$. However, in practice these parabolic schemes are sped up using techniques from~\cite{SchaefferHou}, which allows for potentially much larger time steps to be chosen on the fly and preserves convergence guarantees.



%
%
%
%
%

\section{Computational Results}\label{sec:numerics}
Here we demonstrate the effectiveness of our method with several computational examples. These include reflector design problems involving an omnidirectional source, discontinuous intensity distributions, and intensity distributions supported on sets with complicated geometries.  In each example, we use Algorithm~\ref{alg:reflector} to construct an approximate reflector $\Sigma^h$.

In order to validate our results, we first use the law of reflection~\eqref{eq:mapping2} to perform approximate (forward or inverse) ray-tracing.  We then construct the resulting intensity patterns via approximation of the conservation of energy equation equation~\eqref{eq:conservation} by \[f_1(x_i) \Delta x_i \approx f_2(y_i) \Delta y_i,\]\label{econs} where $\Delta x_i$ and $\Delta y_i$ are the areas of the Voronoi regions containing $x_i$ and $y_i = T(x_i)$ respectively. 

After performing ray tracing, the presence of numerical artifacts may require that the data be post-processed to show the results clearly. This is done by rescaling the colorbars to cut off a very small number of the highest values. Any numerical artifacts are presented in plots of the difference between the desired and ray-traced intensities.

All computations were performed on a 13-inch MacBook Pro, 2.3 GHz Intel Core i5 with 16GB 2133 MHz LPPDDR3 using Matlab R2017b. Each computation utilized around $N\approx 20,000$ points on the sphere.  Where applicable, regularization was performed using $\epsilon=0.3$. The precomputation step of approximating all directional derivatives for $N\approx20,000$ points took about $10$ minutes. Solving the parabolic scheme to find the solution took around $30$ minutes. Ongoing work will develop faster, more accurate versions of this method.  We see therefore that the proposed numerical method can certainly accommodate higher precision computations if necessitated by real-world applications.

\subsection{Peanut Reflector}

\begin{figure}
	\subfigure[Headlight intensity $f_1$ to constant intensity $f_2$]{\includegraphics[width=0.45\textwidth]{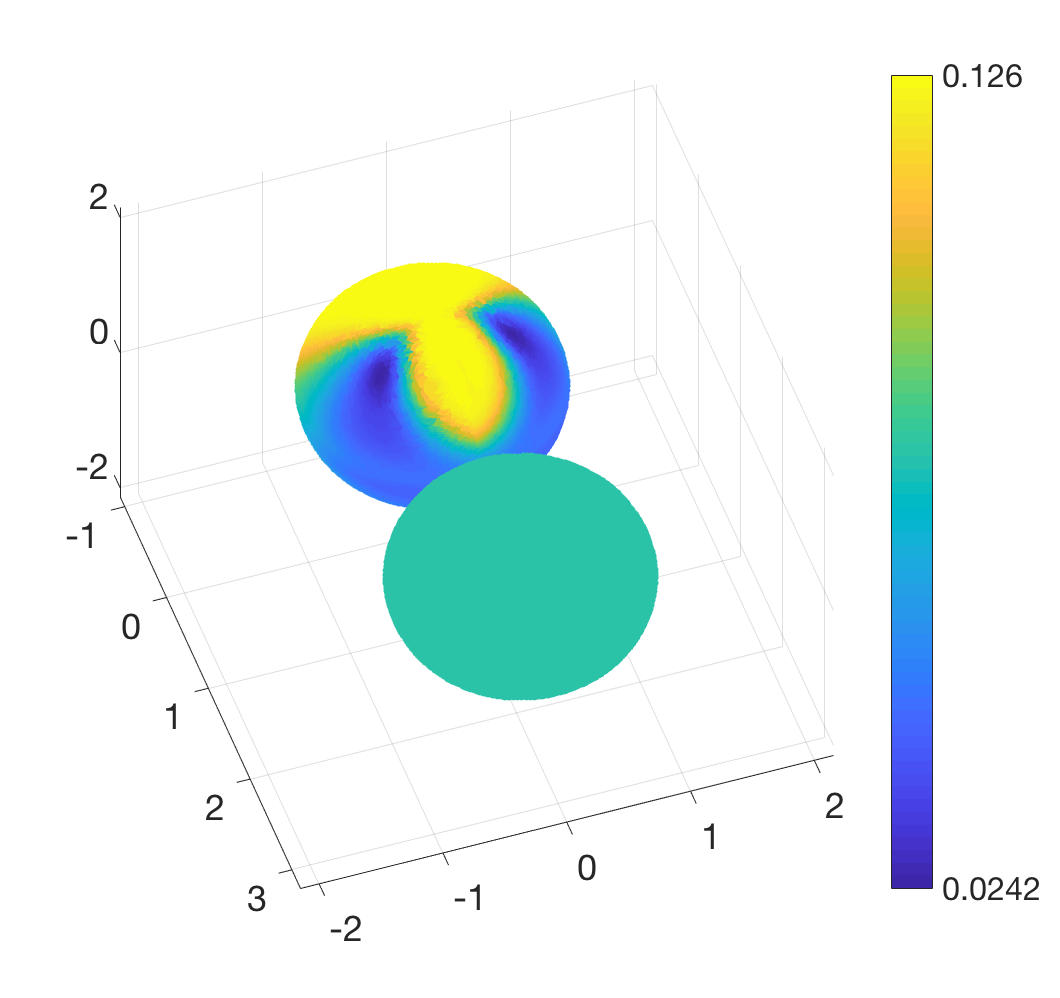}}
		\subfigure[Computed reflector]{\includegraphics[width=0.45\textwidth]{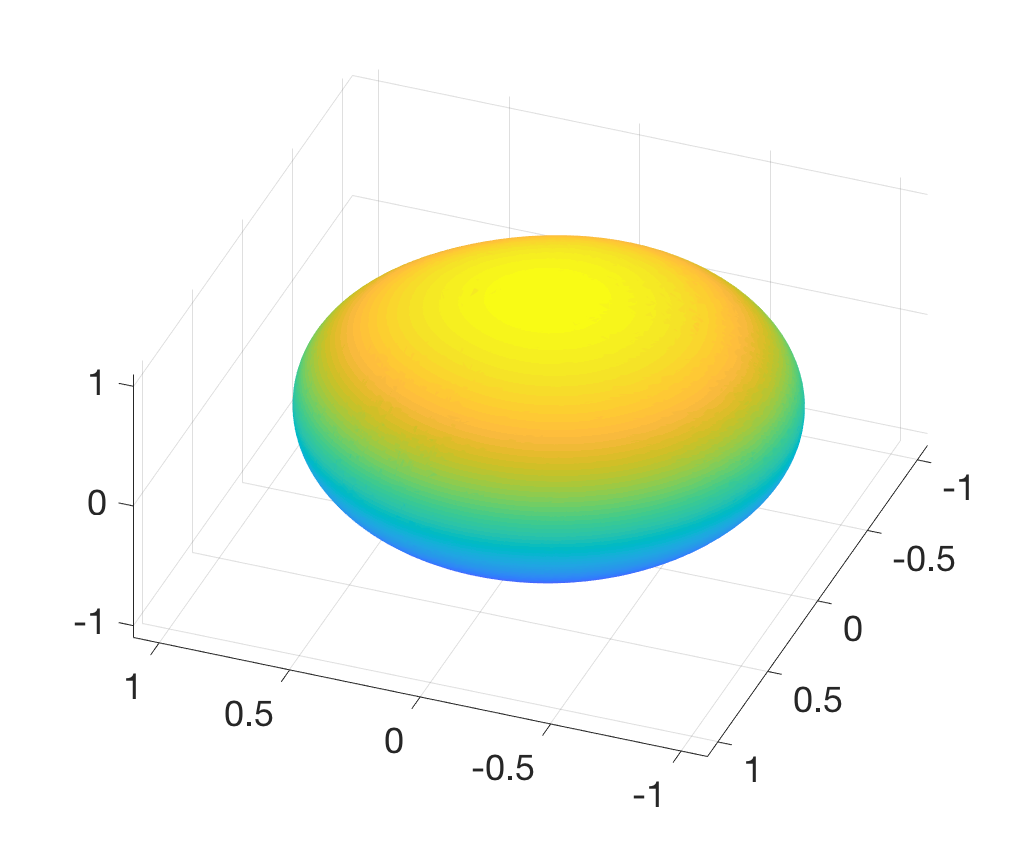}}
	\subfigure[Inverse ray-traced intensity]{\includegraphics[width=0.45\textwidth]{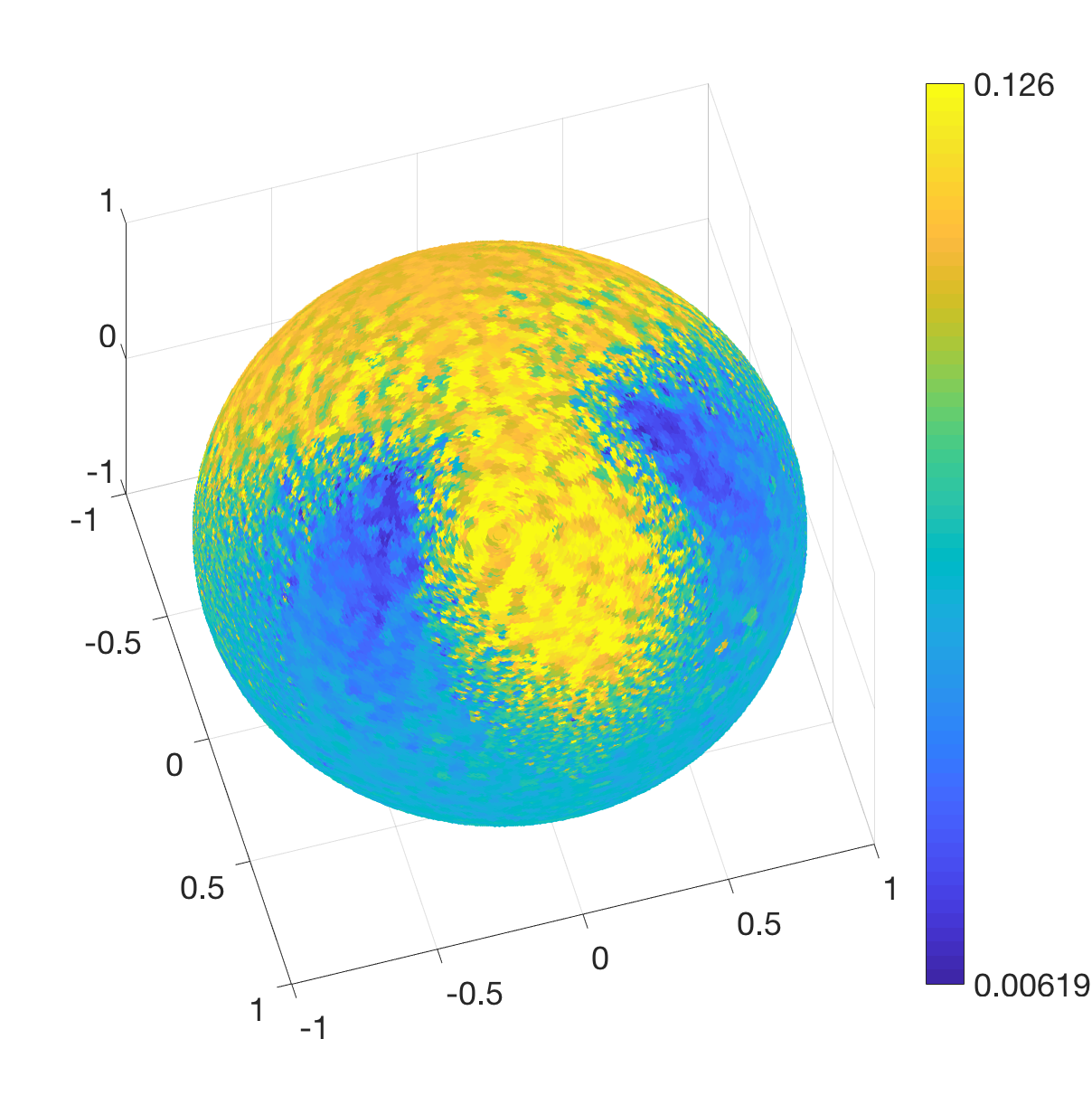}}
	\subfigure[Difference between $f_1$ and inverse ray-traced intensity, with average error of $0.0137$.]{\includegraphics[width=0.45\textwidth]{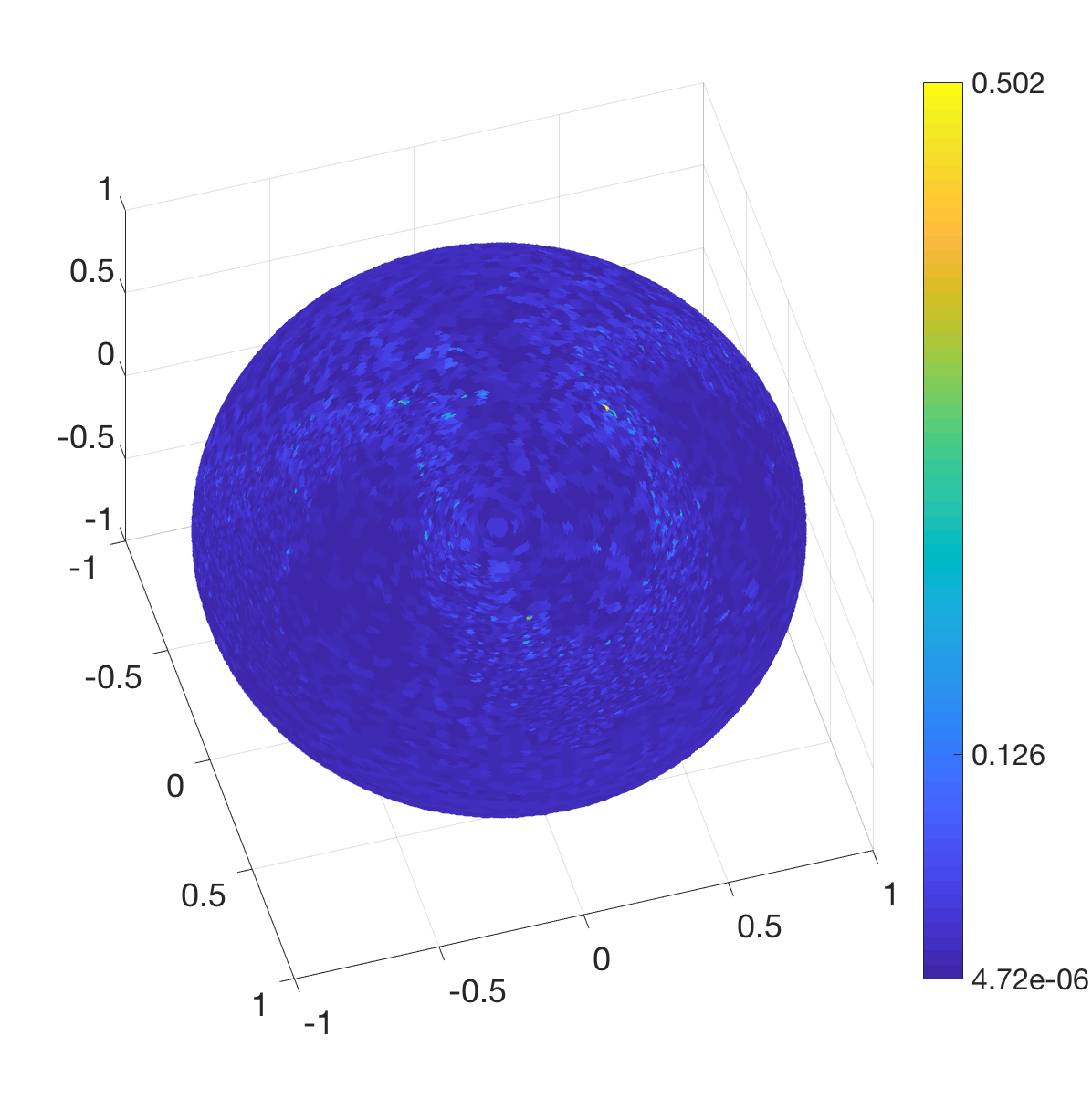}}
	\caption{``Peanut" reflector}\label{fig:peanutresult}
\end{figure}

Following the example of~\cite{RomijnSphere}, we consider a source density coming from an ideal headlight intensity emitting from a vehicle's high beams. This headlight intensity pattern is then mapped to the sphere, and inverted, which becomes the source intensity $f_1$. The target density $f_2$ is constant. 
The computation yields a peanut-shaped oblong reflector lens; see Figure~\ref{fig:peanutresult}. Despite the fact that we anticipate error in the reverse ray trace due to the approximate conservation of energy equation~\eqref{econs}, we see that the absolute error performs quite well in this smooth example.  The average error in the reconstruction is 11\% of the maximum intensity.


\subsection{Discontinuous intensities}
\begin{figure}
	\subfigure[Intensities $f_1$ and $f_2$]{\includegraphics[width=0.45\textwidth]{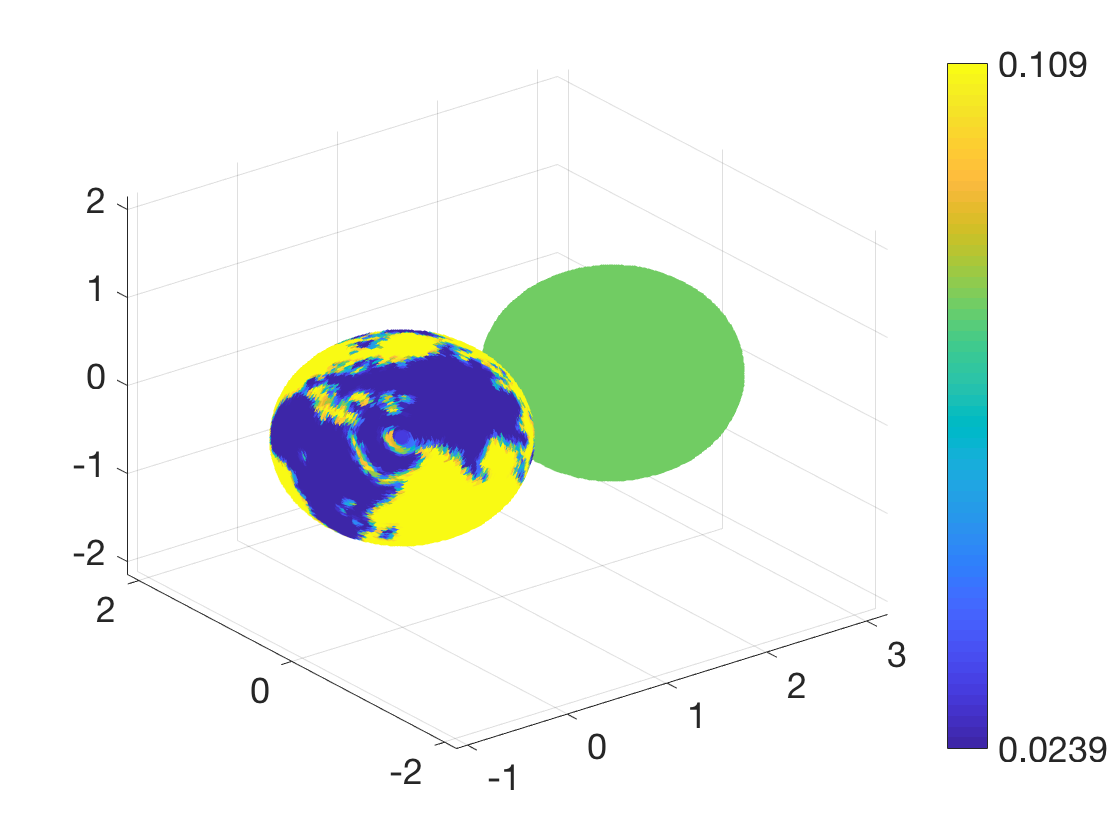}}
	\subfigure[Computed reflector]{\includegraphics[width=0.45\textwidth]{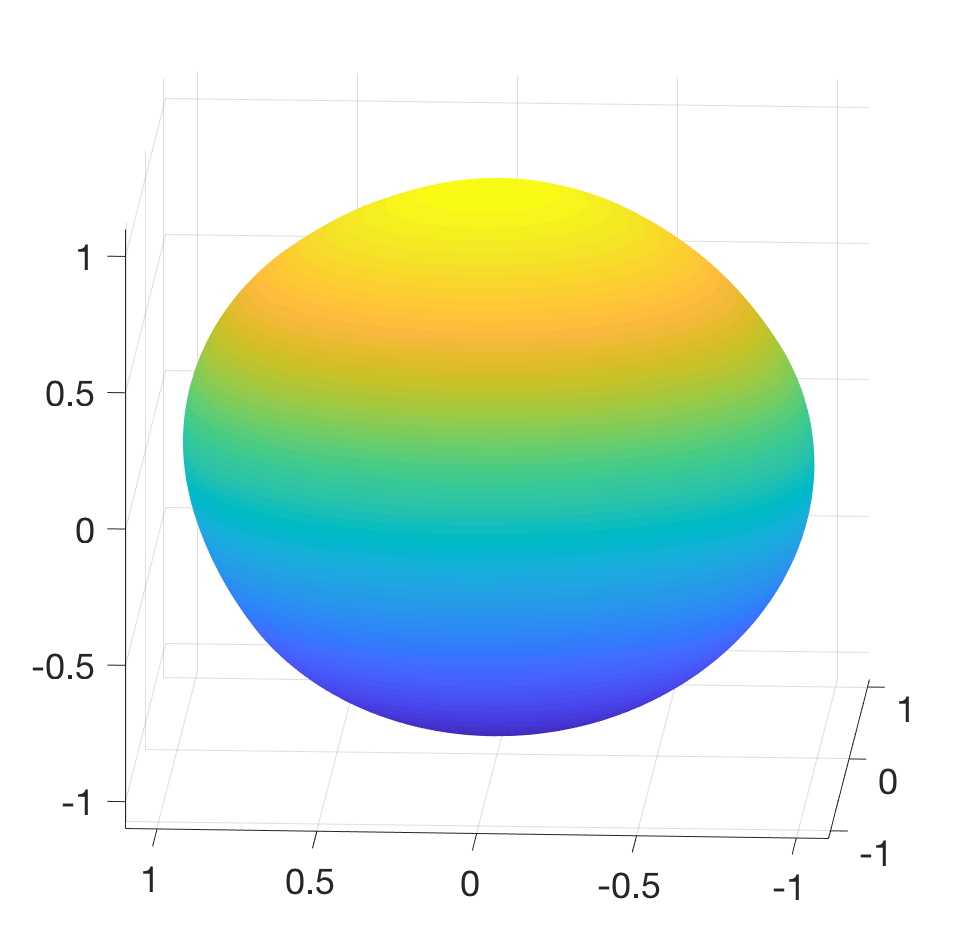}}
	\subfigure[Inverse ray-traced intensity]{\includegraphics[width=0.45\textwidth]{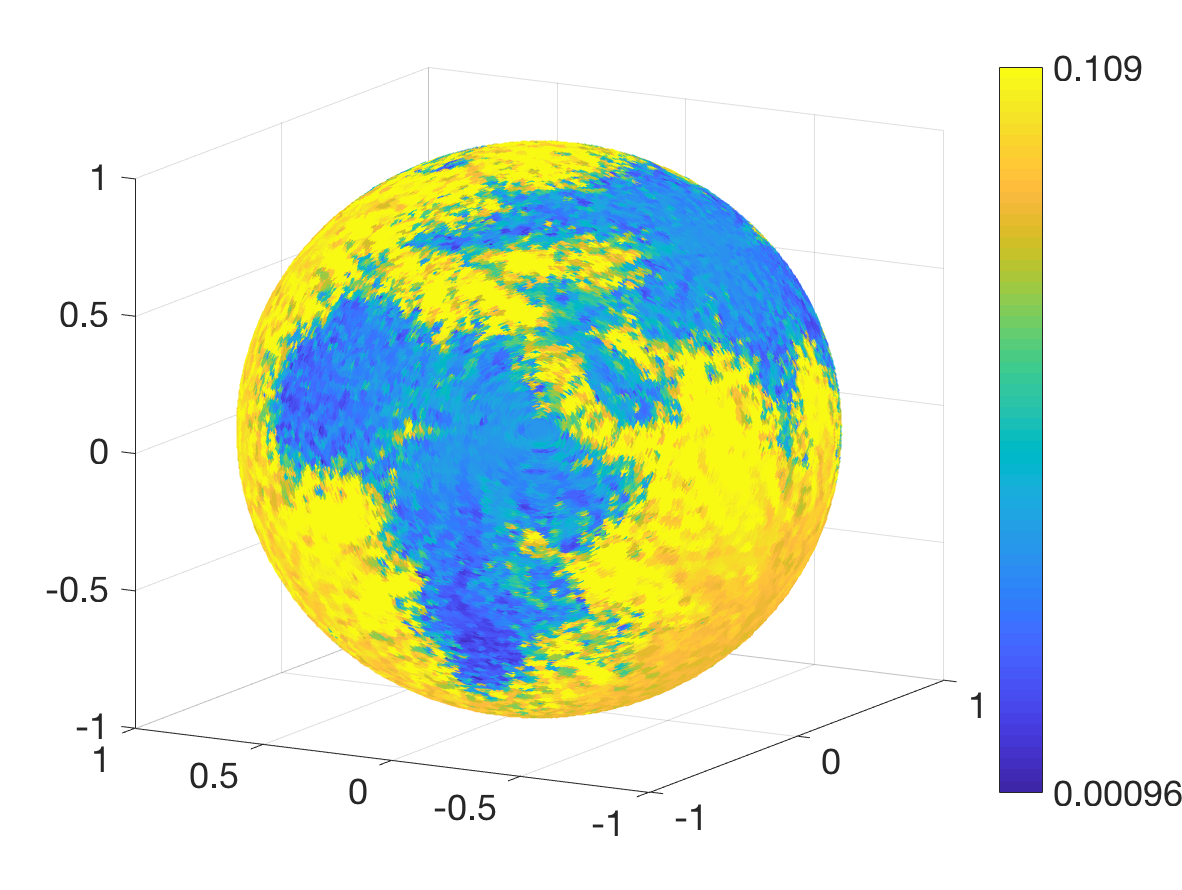}}
	\subfigure[Difference between $f_1$ and inverse ray-traced intensity, with average $L^1$ error of $0.0206$]{\includegraphics[width=0.45\textwidth]{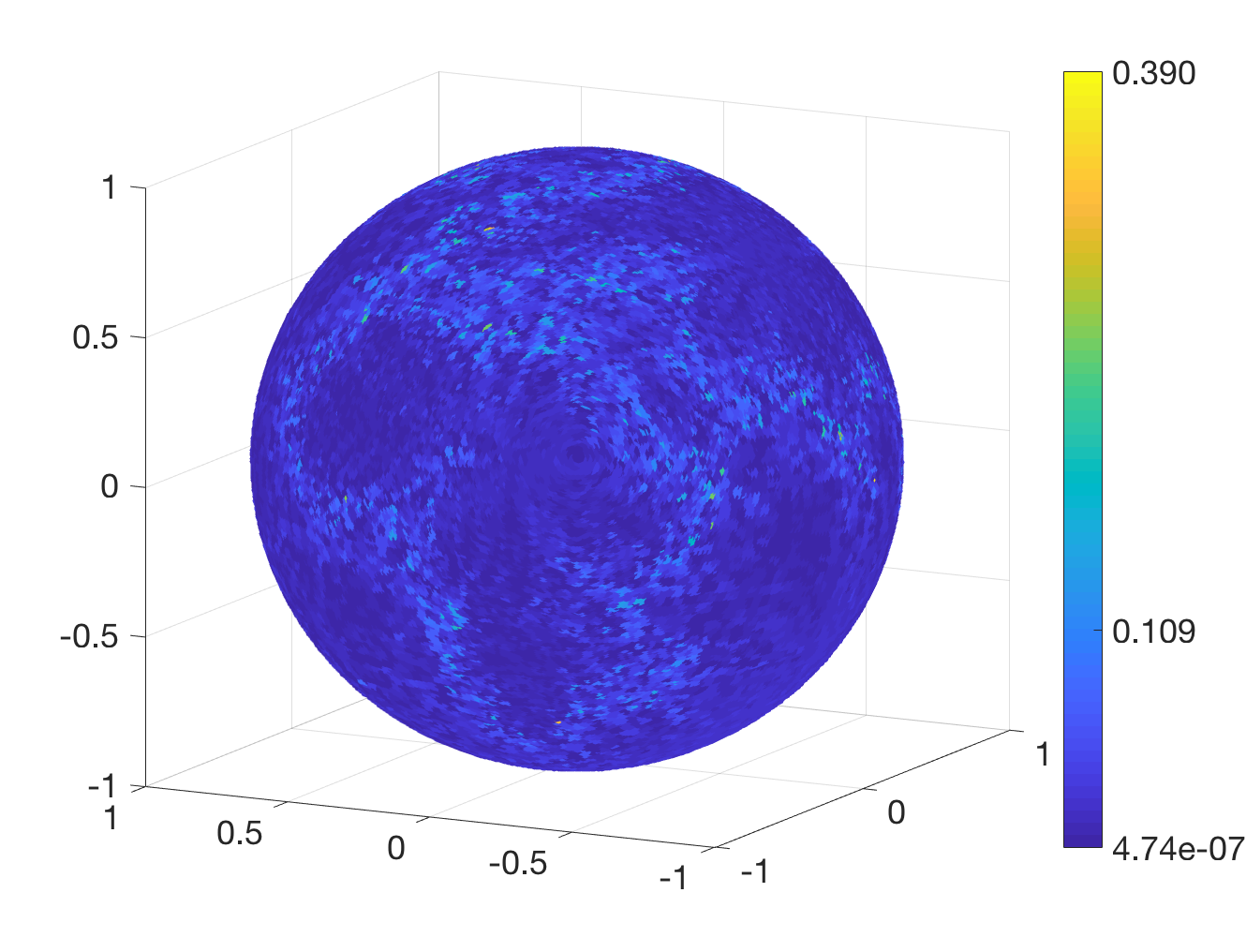}}
	\caption{Discontinuous intensities}\label{fig:globe}
\end{figure}

Next, we demonstrate the effectiveness of our method in dealing with discontinuities and complicated densities. In this example, a discontinuous source mass $f_1$ resembling an inverted map of the world is mapped to a constant density $f_2$; see Figure~\ref{fig:globe}. This is a particularly challenging example given the very complicated structure of the discontinuities.  Nevertheless, we achieve a reconstruction that visually agrees with the world map, with an average error of 19\% of the maximum intensity.

\subsection{Donut intensities}
To further demonstrate the flexibility of our method, we consider the source and target intensities propagating in a donut shape, with a dark region in the center.  These are given by
\begin{equation}
f_1(x,y,z) =
\begin{cases}
\frac{1}{(4\pi/15)(\sqrt{2} + 2)} \left(-4\sqrt{x^2 + y^2} z^3 + 4(x^2 + y^2)^{3/2}z \right), & \sqrt{2}/2 \geq z\geq 0 \\
0, & \text{otherwise}
\end{cases}
\end{equation}
and
\begin{equation}
f_2(x,y,z) =
\begin{cases}
\frac{1}{(4\pi/15)(\sqrt{2} + 2)} \left(-4\sqrt{x^2 + y^2} z^3 + 4(x^2 + y^2)^{3/2}z \right), & 0 \geq z \geq -\sqrt{2}/2 \\
0, & \text{otherwise}
\end{cases}
\end{equation}

\begin{figure}
	\subfigure[Intensities $f_1$ and $f_2$]{\includegraphics[width=0.45\textwidth]{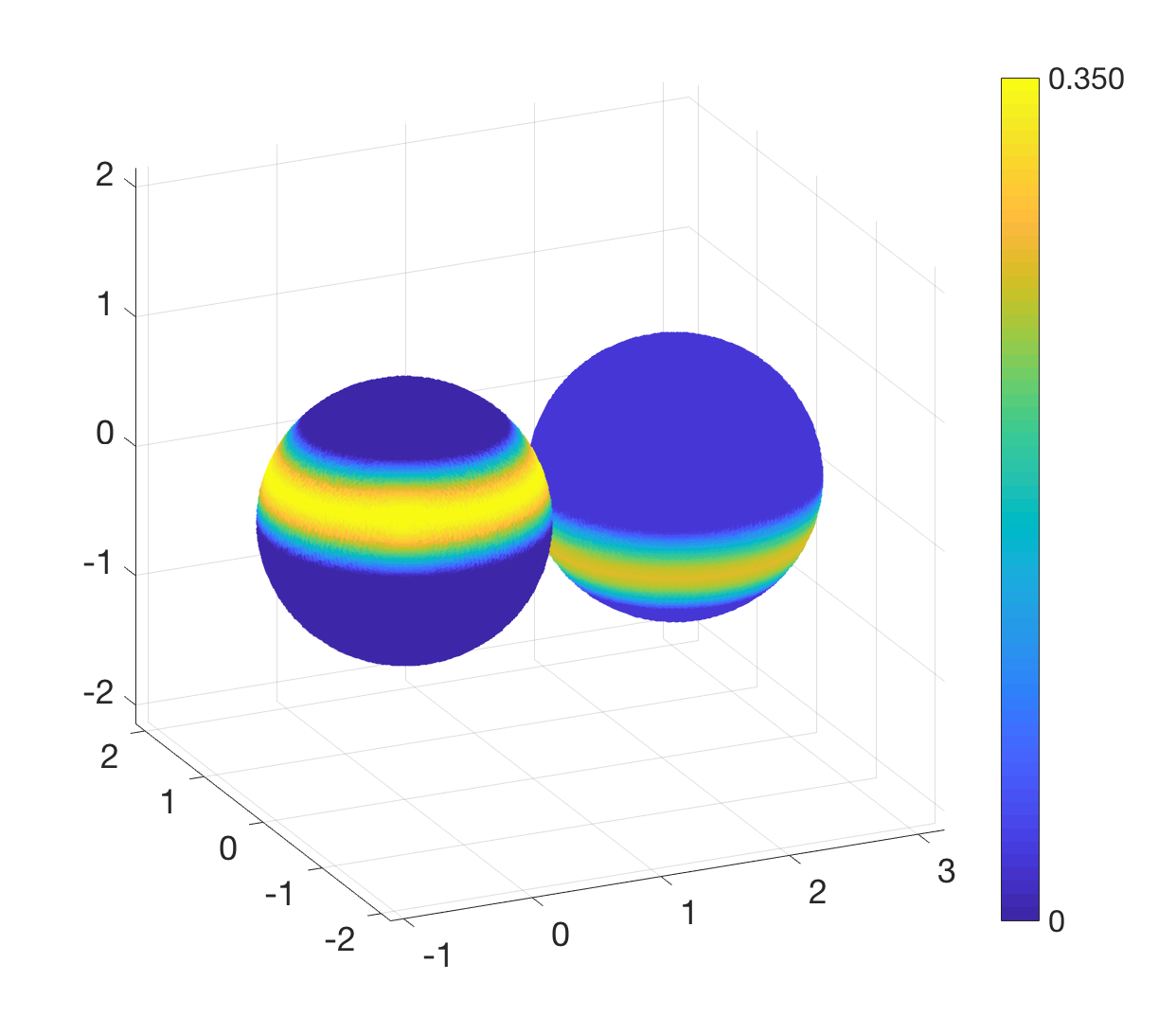}}
	\subfigure[Computed reflector]{\includegraphics[width=0.45\textwidth]{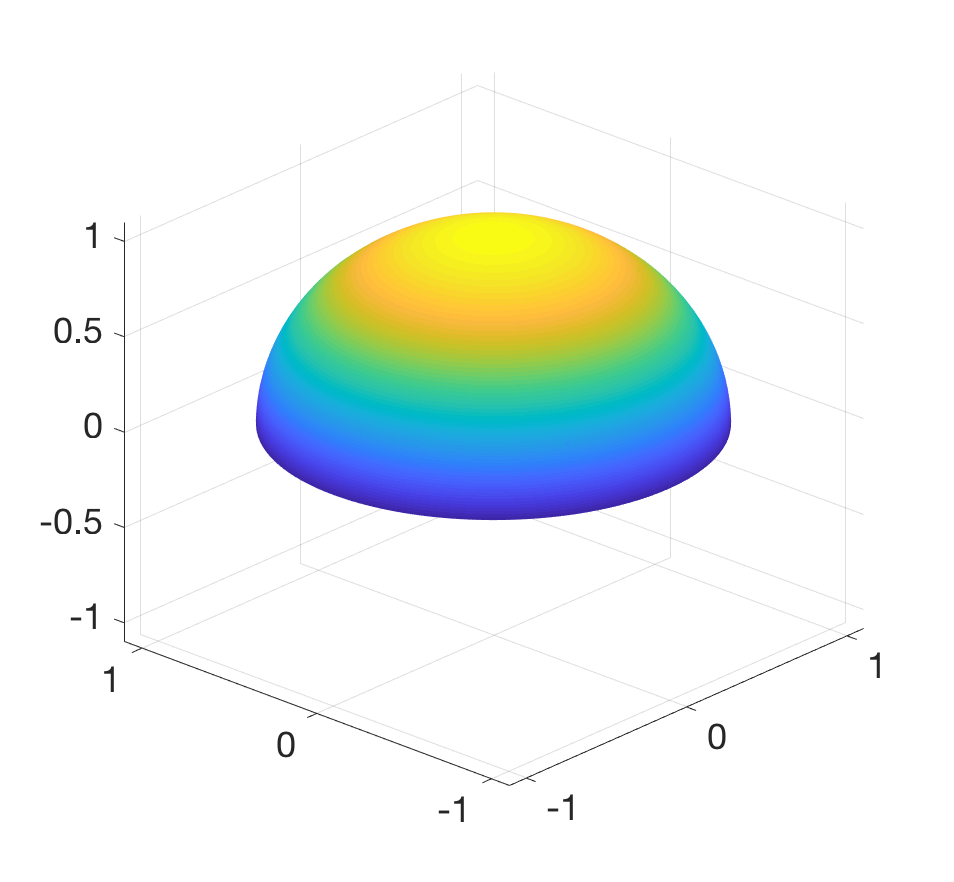}}
	\subfigure[Forward ray-traced intensity]{\includegraphics[width=0.45\textwidth]{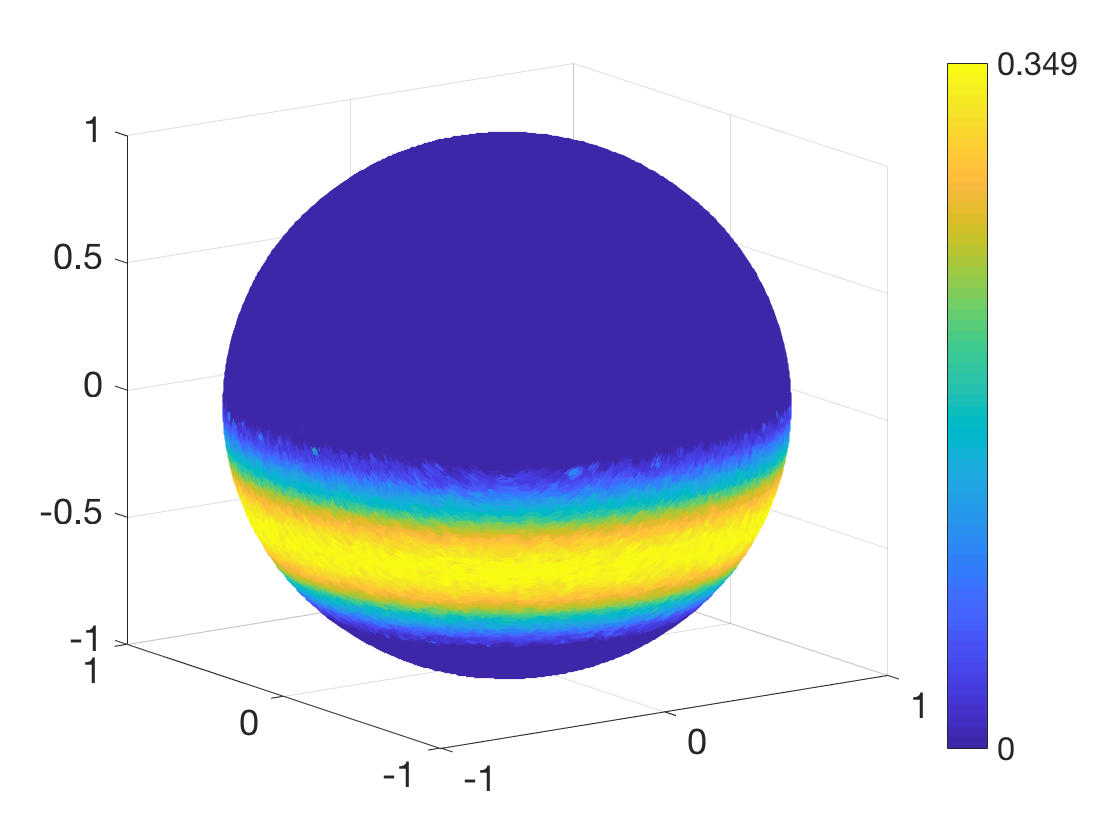}}
	\subfigure[Difference between $f_2$ and forward ray-traced intensity expressed as a percentage of the maximum of $f_2$. Average $L^1$ error of $0.0304$]{\includegraphics[width=0.45\textwidth]{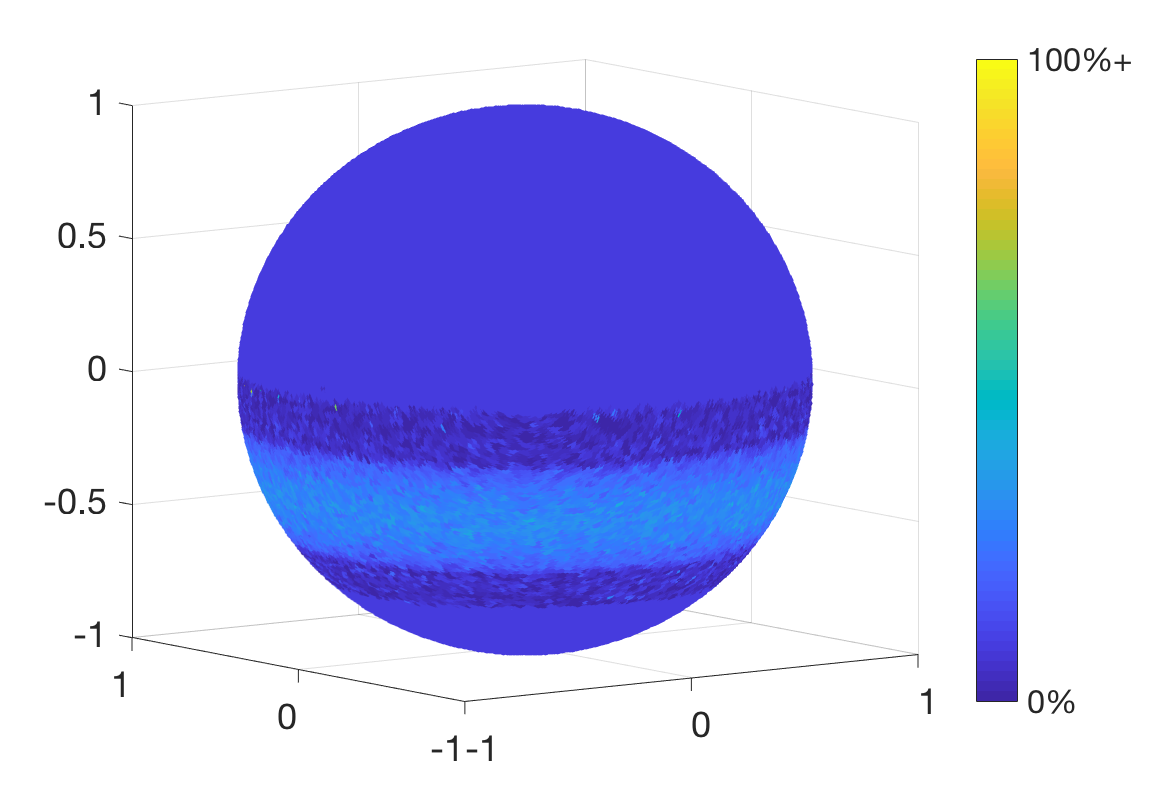}}
	\caption{``Donut" intensities}\label{fig:donutresult}
\end{figure}

These intensities have very complicated support containing holes, which is particularly challenging numerically.  Indeed, this challenge is inherent in the theory of the optimal transport problem.  We note that the $c$-convexity constraint~\eqref{eq:cconvex} requires the domain $\Omega$ to be $c$-convex in order to guarantee construction of the physically relevant solution of the PDE~\eqref{eq:OTPDE}.  Consequently, PDE based methods that are posed only on the support $\Omega$ of the intensity (rather than being extended into the dark regions) will not be assured of producing the correct reflector.  This issue is handled naturally by our method, which is posed on the entire sphere.  Despite the difficulty of this example, our method performs very well, as evidenced in the results of the ray-tracing.  See Figure~\ref{fig:donutresult}.  Average error is 9\% of the maximal intensity.

\subsection{Singular reflector}
We conclude with an example of a hemispheric light source (here designated as $f_2$) that is to be reshaped into a geodesic triangle on the sphere (here designated as $f_1$).  We remark that given the complicated (non $c$-convex) support of this target, we are not even guaranteed the existence of a smooth ($C^1$) reflector; see~\cite{Loeper_OTonSphere}. 

The intensities are defined as follows.  We begin by forming a geodesic triangle $T_{\theta} \subset \Sf$ from the three vertices $(t_{0, \theta}, t_{1, \theta}, t_{2, \theta})$, where we define $t_{j, \theta} = \left( \sin \theta \cos (2\pi j/3), \sin \theta \sin (2\pi j/3), \cos \theta \right)$ for $\pi/2\leq \theta <\pi$. The geodesic triangle is formed by the small region enclosed by the three vertices $t_i$, which are connected by geodesics on the sphere. That is, a point $x_0 \in T_{\theta}$ if $x_0$ satisfies the following three inequalities:
\begin{align*}
x_0 \cdot \left( t_{1,\theta} \times t_{2,\theta} \right) &\leq 0 \\
x_0 \cdot \left( t_{2,\theta} \times t_{3,\theta} \right) &\leq 0 \\
x_0 \cdot \left( t_{3,\theta} \times t_{1,\theta} \right) &\leq 0
\end{align*}

Then the triangular intensity is defined by
\begin{equation}
f_1(x,y,z) = 
\begin{cases}
1/A , \ \ \ &(x,y,z) \in T_{\theta} \\
0, &(x,y,z) \notin T_{\theta}
\end{cases}
\end{equation}
where $A$ is the area of the geodesic triangle $T_{\theta}$ and $\theta = 2.1$.

The second intensity is a smoothed version of the identity function on the northern hemisphere:
\begin{equation}
f_2(x,y,z) = 
\begin{cases}
\frac{2\pi  \log \left( \cosh(a) \right)}{a}\tanh(a z), &z \geq 0 \\
0, & z<0
\end{cases}
\end{equation}
where $a=10$.

For ease of implementation, we perform pre-processing to bound both $f_1$ and $f_2$ away from zero.



\begin{figure}
	\subfigure[Intensities $f_1$ and $f_2$ from below]{\includegraphics[width=0.49\textwidth]{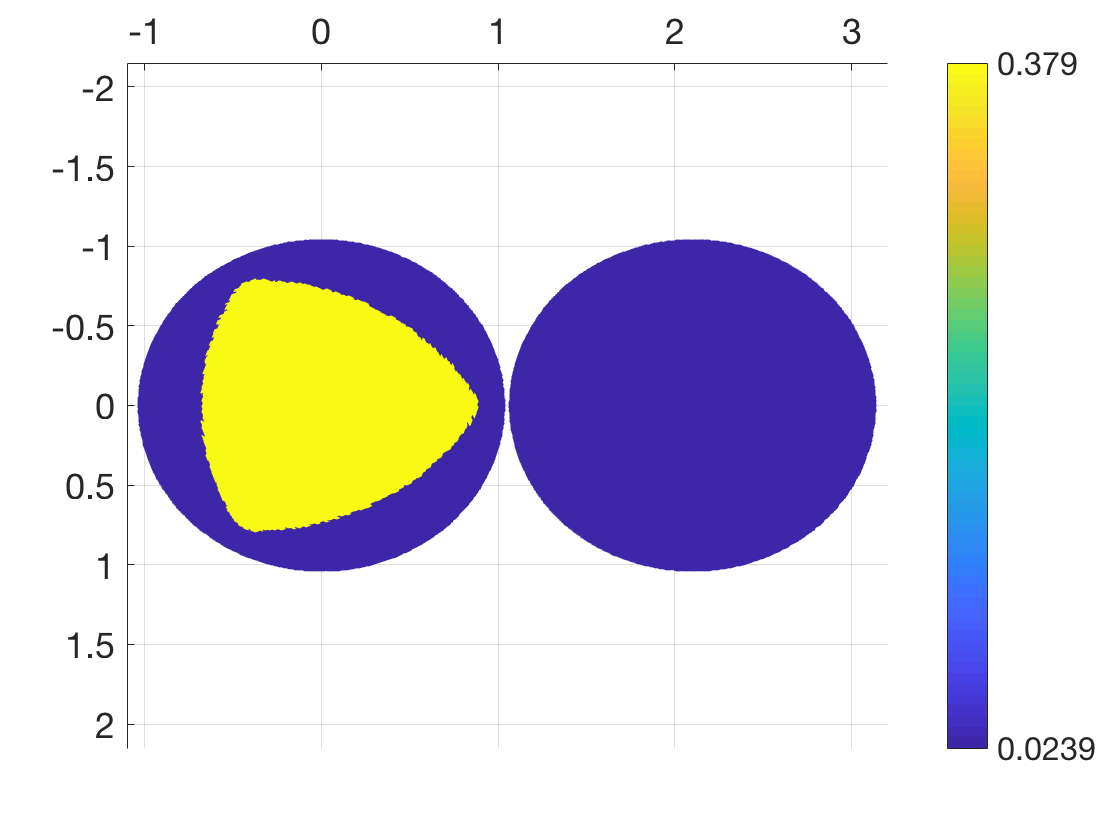}}
	\subfigure[Intensities $f_1$ and $f_2$ from side]{\includegraphics[width=0.49\textwidth]{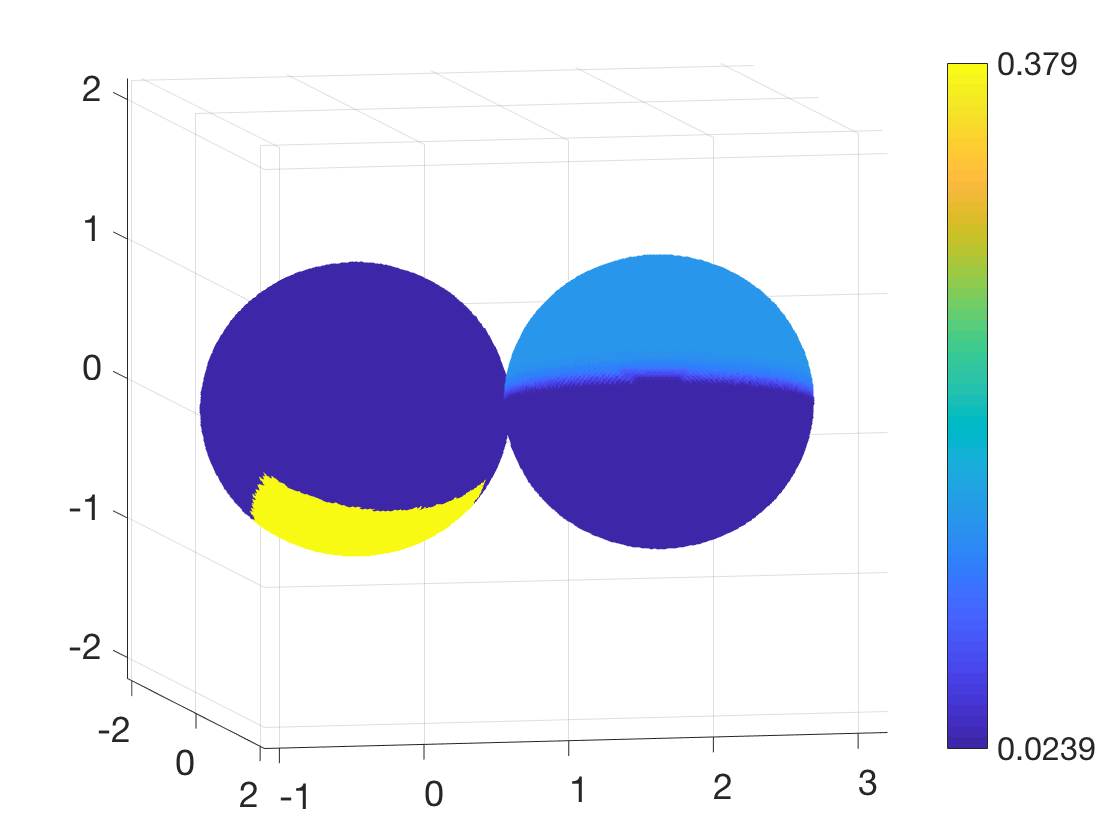}}
	\subfigure[Solution $u^h$]{\includegraphics[width=0.45\textwidth]{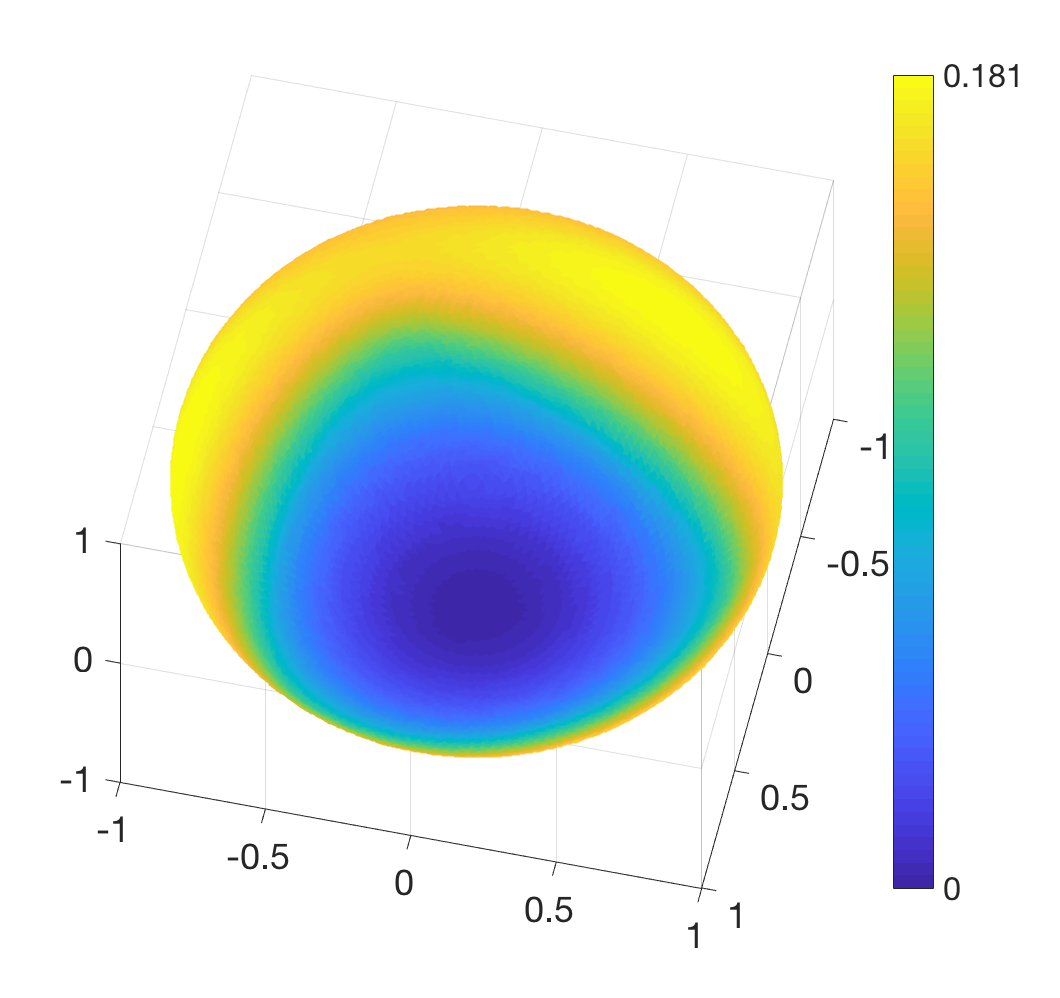}}
	\subfigure[Computed reflector]{\includegraphics[width=0.45\textwidth]{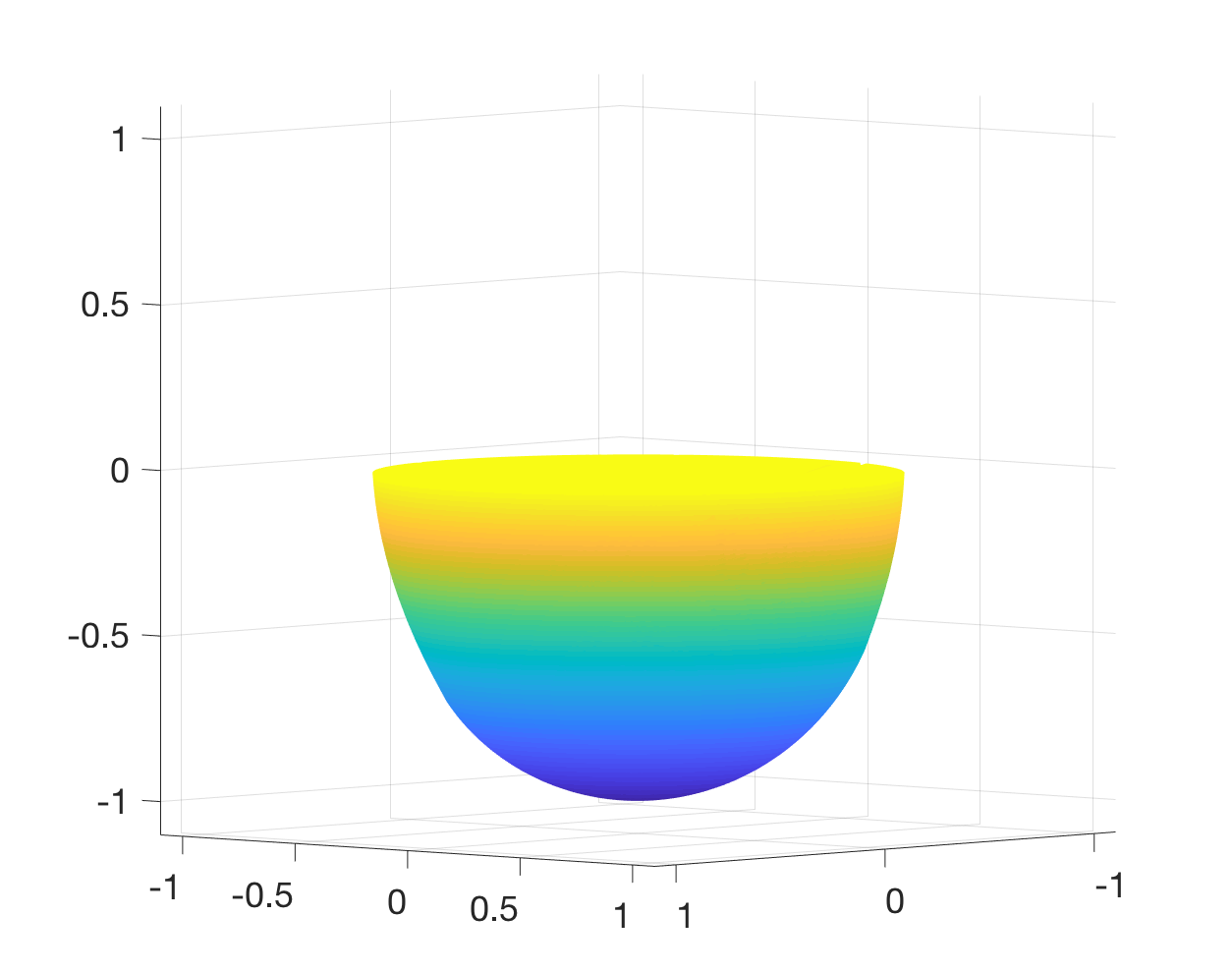}}
	\subfigure[Forward ray-traced intensity]{\includegraphics[width=0.45\textwidth]{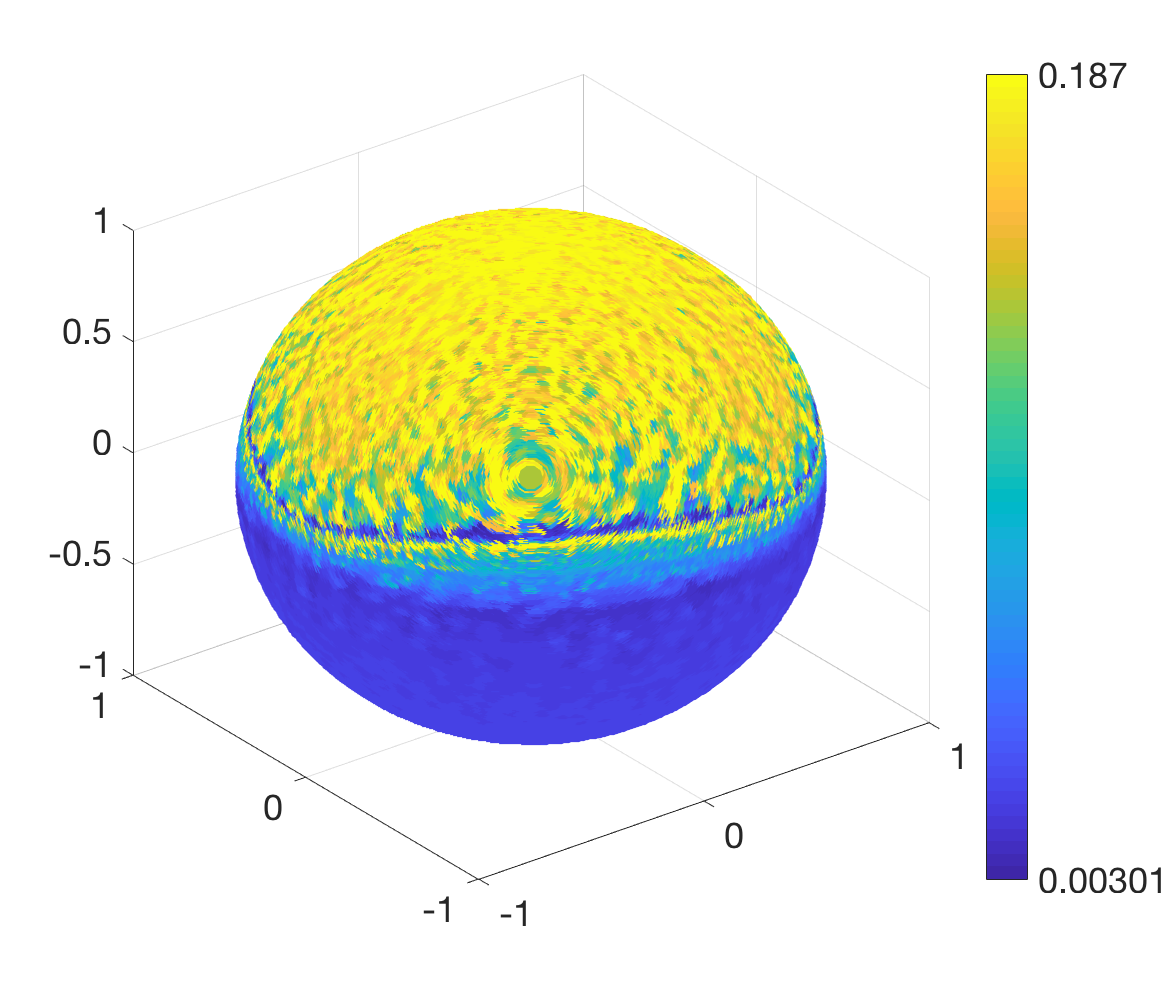}}
	\subfigure[Difference between $f_2$ and forward ray-traced intensity expressed as a percentage of the maximum of $f_2$. Average $L^1$ error of $0.0288$]{\includegraphics[width=0.45\textwidth]{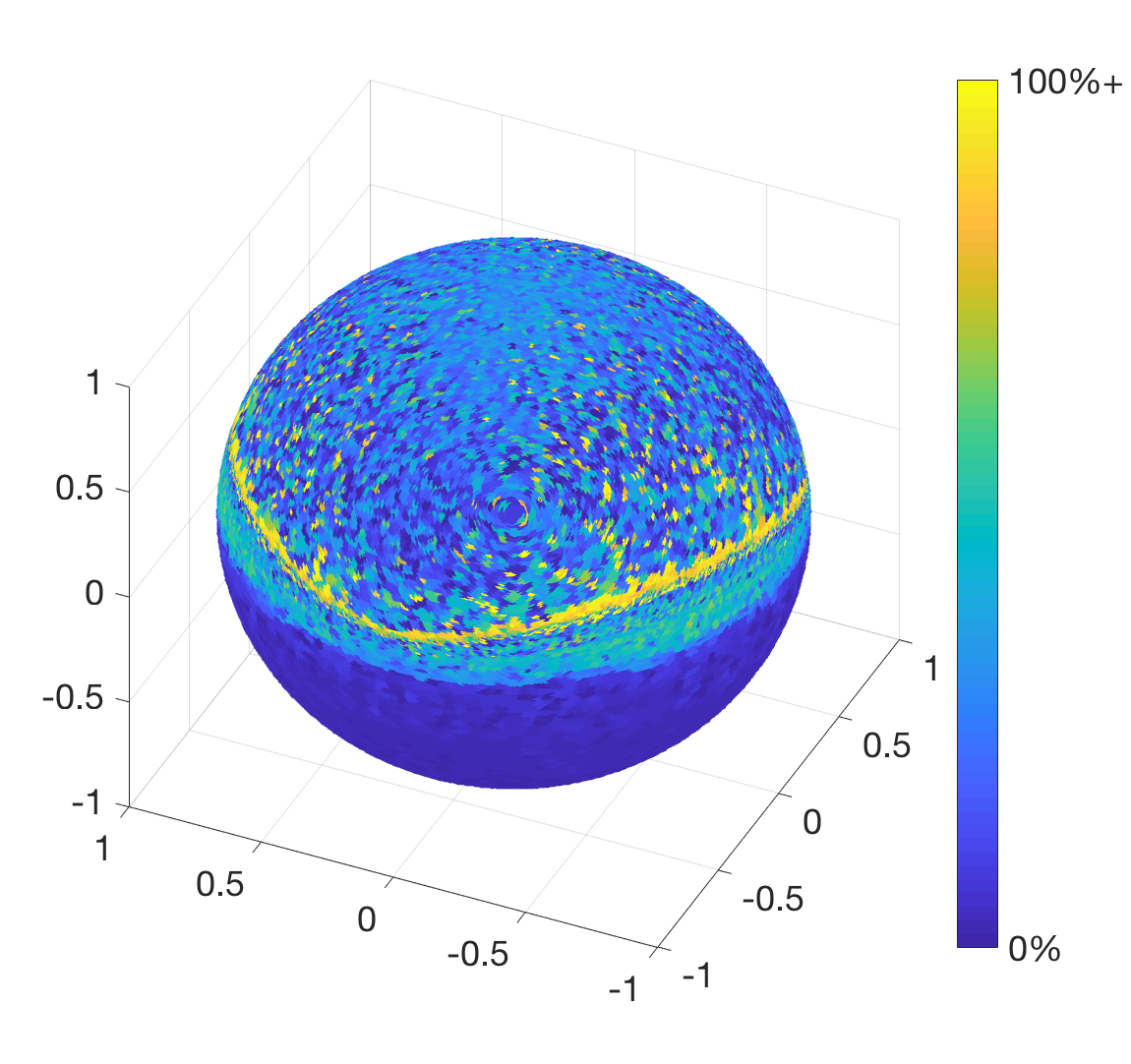}}
	\caption{Singular reflector}\label{fig:triangle}
\end{figure}

Results are presented in Figure~\ref{fig:triangle}.  In the computed reflector, and resulting ray-traced intensity, we observe an approximate triangle shape as expected. In this case, there are notable artifacts present near the boundary of the triangle. However, to some extent these are a limitation of the physics rather than of our method.  We remark that there is no reason to expect the reflector we are approximating to be continuously differentiable, so the accuracy of the ray-tracing verification test is itself rather suspect here.  Nevertheless, the absolute error as compared with the ray trace from the approximate conservation of energy equation mostly performs well, with an average error of 16\% of the maximal intensity.

In a challenging problem like this, where the physics itself may not allow for the existence of a reflector with nice properties (from the perspective of manufacturing and outcome), it may also be useful to view our method as a robust way of obtaining a good approximation of the desired reflector.  This could then be used to initialize an end-game method, not based on optimal transport, that would optimize the reflector surface and enforce any desired smoothness. 

\section{Conclusion}\label{sec:conclusion}
We have introduced a new numerical method for solving the reflector antenna design problem.  The method is based on the reformulation of this design problem as an optimal transport problem on the sphere.  This allows the reflector to be described in terms of the solution to a fully nonlinear elliptic PDE of \MA type, posed on the unit sphere.  We describe a provably convergent finite difference method for solving this PDE, which in turn guarantees that the method will correctly approximate the desired reflector.  The method is robust: convergence guarantees hold even for non-smooth data and reflectors.

We validate this new method through several challenging examples, which include intensities that have complicated discontinuities, that propagate over complicated geometries, or that contain a mix of light and dark regions.  The method performs well even in a final example where the physics does not guarantee the existence of a smooth ($C^1$) reflector.

This new finite difference method provides a rigorous foundation upon which faster and more accurate solvers can be designed.  The idea of pairing slower, more robust approximations (to be used in the most singular regions of the domain) with more traditional high-order methods has been successfully applied to the \MA equation in Euclidean space~\cite{FO_FilteredSchemes}.  In the future, we hope to adapt these techniques to the reflector antenna problem in order to produce higher-quality approximations to the desired reflector surface.

\bibliographystyle{plain}
\bibliography{OTonSphere3}


\end{document}